\documentclass[11pt]{amsart2000}
\usepackage{amscd,amssymb}
\newtheorem{thm}{Theorem}[section]

\newtheorem{prop}[thm]{Proposition}

\theoremstyle{definition}

\theoremstyle{remark}
\newtheorem{remark}[thm]{Remark}

\newtheorem{example}[thm]{Example}

\newtheorem{Question}[thm]{Question}

\numberwithin{equation}{section}
\def\R {{\mathbb R}}
\def\C {{\mathbb C}}
\def\T {{\mathbb T}}
\def\Z {{\mathbb Z}}
\def\e{{\varepsilon}}
\def\Iso{{\mathrm{Is}}\,}
\def\diam{{\mathrm{diam}}}
\def\Exp{\operatorname{Exp}}
\def\Symm{\operatorname{Symm}}
\def\cl#1{\overline {#1}}
\def\sbs{\subset}
\def\stm{\setminus}
\def\k{\kappa}
\def\G{\Gamma}
\def\o{\omega}
\def\sB{\mathcal{B}}
\def\sL{\mathcal{L}}
\def\sR{\mathcal{R}}
\def\sS{\mathcal{S}}
\def\sU{\mathcal{U}}
\def\sN{\mathcal{N}}
\def\RUC{{\hbox{RUC\,}^b}}
\def\obr{^{-1}}

\begin{document}
\setlength{\unitlength}{0.01in}
\linethickness{0.01in}
\begin{center}
\begin{picture}(474,66)(0,0) 
\multiput(0,66)(1,0){40}{\line(0,-1){24}}
\multiput(43,65)(1,-1){24}{\line(0,-1){40}}
\multiput(1,39)(1,-1){40}{\line(1,0){24}}
\multiput(70,2)(1,1){24}{\line(0,1){40}}
\multiput(72,0)(1,1){24}{\line(1,0){40}}
\multiput(97,66)(1,0){40}{\line(0,-1){40}} 
\put(143,66){\makebox(0,0)[tl]{\footnotesize Proceedings of the Ninth Prague Topological Symposium}}
\put(143,50){\makebox(0,0)[tl]{\footnotesize Contributed papers from the symposium held in}}
\put(143,34){\makebox(0,0)[tl]{\footnotesize Prague, Czech Republic, August 19--25, 2001}}
\end{picture}
\end{center}
\vspace{0.25in}
\setcounter{page}{331}
\title{Compactifications of topological groups}
\author{Vladimir Uspenskij}
\address{Department of Mathematics,
321 Morton Hall, Ohio University, Athens, Ohio 45701, USA}
\email{uspensk@bing.math.ohiou.edu}
\thanks{The author was an invited speaker at the Ninth Prague Topological Symposium.}
\thanks{Vladimir Uspenskij,
{\em Compactifications of topological groups},
Proceedings of the Ninth Prague Topological Symposium, (Prague, 2001),
pp.~331--346, Topology Atlas, Toronto, 2002}
\begin{abstract} 
Every topological group $G$ has some natural compactifications which can 
be a useful tool of studying $G$. 
We discuss the following constructions:
(1) the {\it greatest ambit} $S(G)$ is the compactification corresponding
to the algebra of all right uniformly continuous bounded functions on $G$;
(2) the {\it Roelcke compactification\/} $R(G)$ corresponds to the algebra
of functions which are both left and right uniformly continuous; 
(3) the {\it weakly almost periodic\/} compactification $W(G)$ is the 
envelopping compact semitopological semigroup of $G$ (`semitopological'
means that the multiplication is separately continuous). 
The {\it universal minimal compact $G$-space} $X=M_G$ is characterized by
the following properties: 
(1) $X$ has no proper closed $G$-invariant subsets; 
(2) for every compact $G$-space $Y$ there exists a $G$-map $X\to Y$. A
group $G$ is {\it extremely amenable}, or has the {\it fixed point on
compacta} property, if $M_G$ is a singleton. 
We discuss some results and questions by V.~Pestov and E.~Glasner on
extremely amenable groups.

The Roelcke compactifications were used by M.~Megrelishvili to prove that
$W(G)$ can be a singleton. They can be used to prove that certain groups 
are minimal. A topological group is {\em minimal} if it does not admit a
strictly coarser Hausdorff group topology. 
\end{abstract}
\maketitle

\setcounter{tocdepth}{1}

\section{Introduction}\label{secintro}

This is a write-up of the lecture that I gave in the 9th Prague 
Topological Symposium on 24 August 2001.

Every topological group $G$ has some natural compactifications. They can
be described as the maximal ideal spaces of certain function algebras, or
as the Samuel compactifications for certain uniformities on $G$. Some
compactifications of $G$ carry an algebraic structure, and may be useful
for studying the group $G$ itself.

We consider, in particular, the following constructions: the greatest
ambit $S(G)$ and the universal minimal compact $G$-space $M_G$
(Sections~\ref{secsg} and~\ref{secmg}); the Roelcke compactification
$R(G)$ (Section~\ref{secroel}); the weakly almost periodic
compactification $W(G)$ (Section~\ref{secwap}). In the last case the
canonical map $G\to W(G)$ need not be an embedding, as Megrelishvili has
recently proved \cite{MeForum}. In Section~\ref{secury} we discuss the
group of isometries of the Urysohn universal metric space $U$.

There are two important classes of topological groups: the groups of the
form $H(K)$, $K$ compact, and of the form $\Iso(M)$, $M$ metric. For a
compact space $K$ let $H(K)$ be the topological group of all
self-homeomorphism of $K$, with the compact-open topology (which is the
same as the topology of uniform convergence, if $K$ is equipped with its
unique compatible uniformity). For a metric space $M$ let $\Iso(M)$ be the
topological group of all isometries of $M$ onto itself. The topology that
we consider on $\Iso(M)$ is the topology of pointwise convergence or,
equivalently, the compact-open topology -- the two topologies coincide on
$\Iso(M)$. If $B$ is a Banach space, we use the symbol $\Iso_0(B)$ to
denote the topological group of all linear isometries of $B$. The
pointwise convergence topology on $\Iso_0(B)$ is also called the {\em
strong operator topology}. The {\em weak operator topology} on $\Iso_0(B)$
is the topology inherited from the product $(B_w)^B$, where $B_w$ is $B$
with the weak topology. The weak topology on $\Iso_0(B)$ in general is
not compatible with the group structure \cite{MeElena}. On the other
hand, if $B$ is reflexive, the weak and strong topologies on $\Iso(B)$
agree \cite{MeElena}, see Section~\ref{secwap} below.

Every topological group $G$ has a {\it topologically faithful} 
representation by homeomorphisms of a compact space or by isometries of a
Banach space. This means that $G$ is isomorphic (as a topological
group) to a subgroup of $H(K)$ for some compact $K$ and also to a subgroup
of $\Iso_0(B)$ for some Banach space $B$.
(The term `topologically faithful' was suggested in \cite{P3}.)
To see this, take for $B$ the Banach space $\RUC(G)$ of all right 
uniformly continuous bounded complex functions on $G$. 
Recall that a function $f$ on $G$ is {\it right uniformly continuous} if
$$
\forall\e>0\,\,\exists V\in\sN(G)\,\,\forall x,y\in G\,\,
(xy\obr\in V\implies|f(y)-f(x)|<\e),
$$
where $\sN(G)$ is the filter of neighbourhoods of unity. 
The embedding $G\to \Iso_0(B)$ comes from the action of $G$ on $B$ defined
by $gf(x)=f(g\obr x)$ ($g,x\in G$, $f\in B$). 
If $K$ is the unit ball of the dual space $B^*$ with the 
weak$^*$-topology, then the natural map $\Iso_0(B)\to H(K)$ is an
isomorphic embedding, and we obtain a topologically faithful
representation of $G$ by homeomorphisms of $K$.

All {\em maps} are assumed to be continuous, and `compact' includes 
`Hausdorff'. 

Most of the results and ideas presented in this note can be found in the 
excellent survey \cite{P3}.

\section{Greatest ambit $S(G)$}\label{secsg}

Let $G$ be a topological group. The Banach space $B=\RUC(G)$ of all right
uniformly continuous bounded complex functions on $G$ is a $C^*$-algebra,
and $G$ acts on $B$ by $C^*$-algebra automorphisms. Let $S(G)$ be the
(compact) maximal ideal space of $B$. It is the least compactification of
$G$ over which all functions from $B$ can be extended. The topological
group of all $C^*$-algebra automorphisms of $B$ is naturally isomorphic to
$H(S(G))$. It follows that $G$ acts on $S(G)$, and the natural
homomorphism $G\to H(S(G))$ is a topological embedding.

The space $S(G)$ can also be described as the Samuel compactification of
the uniform space $(G, \sR)$. Here $\sR$ is the right uniformity on $G$.
The basic entourages for $\sR$ are of the form $\{(x,y)\in G\times
G:xy\obr\in V\}$, where $V\in \sN(G)$. The {\it Samuel compactification}
of a uniform space $(X,\sU)$ is the completion of $X$ with respect to the
finest precompact uniformity which is coarser than $\sU$.

We shall consider $G$ as a dense subpace of $S(G)$. The action $G\times
S(G)\to S(G)$ extends the multiplication $G\times G\to G$.

A {\em $G$-space} is a topological space $X$ with a continuous action of
$G$, that is, a map $G\times X\to X$ satisfying $g(hx)=(gh)x$ and $1x=x$
($g,h\in G$, $x\in X$). A {\it $G$-map} is a map $f:X\to Y$ between
$G$-spaces such that $f(gx)=gf(x)$ for all $x\in X$, $g\in G$. The
$G$-space $S(G)$ has a distinguished point $e$ (the unity), and the pair
$(S(G), e)$ has the following universal property: for every compact
$G$-space $X$ and every $p\in X$ there exists a unique $G$-map $f:S(G)\to
X$ such that $f(e)=p$. Indeed, the map $g\mapsto gp$ from $G$ to $X$ is
$\sR$-uniformly continuous and hence can be extended over $S(G)$.

The space $S(G)$ (or the pair $(S(G),e)$) is called the {\em greatest
ambit} of $G$. Let us show that $S(G)$ has a natural structure of a
left-topological semigroup. A {\em semigroup} is a set with an associative
multiplication. A semigroup $X$ is {\em left-topological} if it is a
topological space and for every $y\in X$ the self-map $x\mapsto xy$ of $X$
is continuous. (Some authors use the term {\em right-topological} for
this.)

\begin{thm}
For every topological group $G$ the greatest ambit $X=\sS(G)$ has a
natural structure of a left-topological semigroup with a unity such that
the multiplication $X\times X\to X$ extends the action $G\times X\to X$.
\end{thm}

\begin{proof}
Let $x, y\in X$.
In virtue of the universal property of $X$, there is a unique $G$-map
$r_y:X\to X$ such that $r_y(e)=y$. 
Define $xy=r_y(x)$.
Let us verify that the multiplication $(x,y)\mapsto xy$ has the required
properties.
For a fixed $y$, the map $x\mapsto xy$ is equal to $r_y$ and hence is
continuous. 
If $y,z\in X$, the self-maps $r_zr_y$ and $r_{yz}$ of $X$ are equal, since
both are $G$-maps sending $e$ to $yz=r_z(y)$. 
This means that the multiplication on $X$ is associative. The
distinguished element $e\in X$ is the unity of $X$: we have $ex=r_x(e)=x$
and $xe=r_e(x)=x$. 
If $g\in G$ and $x\in X$, the expression $gx$ can be understood in two
ways: in the sense of the exterior action of $G$ on $X$ and as a product
in $X$. 
To see that these two meanings agree, note that 
$r_x(g)=r_x(ge)=gr_x(e)=gx$ 
(the exterior action is meant in the last two terms; the middle equality
holds since $r_x$ is a $G$-map). 
\end{proof}

\section{Universal minimal compact $G$-space}\label{secmg}

Let us define the {\em universal minimal compact $G$-space} $M_G$. A
$G$-space $X$ is {\em minimal} if it has no proper $G$-invariant closed
subsets or, equivalently, if the orbit $Gx$ is dense in $X$ for every
$x\in X$. The universal minimal compact $G$-space $M_G$ is characterized
by the following property: $M_G$ is a minimal compact $G$-space, and for
every compact minimal $G$-space $X$ there exists a $G$-map of $M_G$ onto
$X$. Since Zorn's lemma implies that every compact $G$-space has a minimal
compact $G$-subspace, it follows that for every compact $G$-space $X$,
minimal or not, there exists a $G$-map of $M_G$ to $X$.

The existence of $M_G$ is easy: take for $M_G$ any minimal closed
$G$-subspace of $S(G)$. The universal property of $(S(G), e)$ implies the
corresponding universal property of $M_G$. It is also true that $M_G$ is
unique, in the sense that any two universal minimal compact $G$-spaces are
isomorphic \cite{Aus}. For the reader's convenience we give a proof of
this fact.

Let $X=S(G)$. For $a\in X$ let $r_a$ be the map $x\mapsto xa$ of $X$ to
itself.

\begin{prop}
If $f:X\to X$ is a $G$-self-map and $a=f(e)$, then $f=r_a$.
\label{pr4}
\end{prop}

\begin{proof}
We have $f(x)=f(xe)=xf(e)=xa=r_a(x)$ for all $x\in G$ and hence for all
$x\in X$. 
\end{proof}

A subset $I\sbs X$ is a {\em left ideal} if $XI\sbs I$. Closed
$G$-subspaces of $X$ are the same as closed left ideals of $X$. An element
$x$ of a semigroup is an {\it idempotent\/} if $x^2=x$. Every closed
$G$-subspace of $X$, being a left ideal, is moreover a left-topological
compact semigroup and hence contains an idempotent, according to the
following fundamental result of R.~Ellis (see \cite[Proposition~2.1]{Rup}
or \cite[Theorem~3.11]{BJM}):

\begin{thm}
Every non-empty compact left-topological semigroup $K$ contains an idempotent.
\label{idemp}
\end{thm}

\begin{proof}
Zorn's lemma implies that there exists a minimal element $Y$ in the set
of all closed non-empty subsemigroups of $K$. Fix $a\in Y$. We claim that
$a^2=a$ (and hence $Y$ is a singleton). The set $Y\!a$, being a closed
subsemigroup of $Y$, is equal to $Y$. It follows that the closed
subsemigroup $Z=\{x\in Y: xa=a\}$ is non-empty. Hence $Z=Y$ and $xa=a$ for
every $x\in Y$. In particular, $a^2=a$.
\end{proof}

Let $M$ be a minimal closed left ideal of $X$. We have just proved that
there is an idempotent $p\in M$. Since $Xp$ is a closed left ideal
contained in $M$, we have $Xp=M$. It follows that $xp=x$ for every $x\in
M$. The $G$-map $r_p:X\to M$ defined by $r_p(x)=xp$ is a retraction of $X$
onto $M$.

\begin{prop}
Every $G$-map $f:M\to M$ has the form $f(x)=xy$ for some $y\in M$.
\label{pr2}
\end{prop}

\begin{proof}
The composition $h=fr_p:X\to M$ is a $G$-map of $X$ into itself, hence it
has the form $h=r_y$, where $y=h(e)\in M$ (Proposition~\ref{pr4}). 
Since $r_p\restriction M =\text{Id}$, we have 
$f=h\restriction M=r_y\restriction M$.
\end{proof}

\begin{prop}
Every $G$-map $f:M\to M$ is bijective.
\label{pr3}
\end{prop}

\begin{proof}
According to Proposition~\ref{pr2}, there is $a\in M$ such that $f(x)=xa$ 
for all $x\in M$. 
Since $Ma$ is a closed left ideal of $X$ contained in $M$, we have $Ma=M$
by the minimality of $M$. 
Thus there exists $b\in M$ such that $ba=p$. 
Let $g:M\to M$ be the $G$-map defined by $g(x)=xb$. 
Then $fg(x)=xba=xp=x$ for every $x\in M$, therefore $fg=1$ (the identity
map of $M$). 
We have proved that in the semigroup $S$ of all $G$-self-maps of $M$,
every element has a right inverse. 
Hence $S$ is a group. 
(Alternatively, we first deduce from the equality $fg=1$ that all elements
of $S$ are surjective and then, applying this to $g$, we see that $f$ is
also injective.)
\end{proof}

We are now in a position to prove

\begin{thm}
Every universal compact minimal $G$-space is isomorphic to $M$.
\label{uniqmin}
\end{thm}

\begin{proof}
We noted that the minimal compact $G$-space $M$ is itself universal: if
$Y$ is any compact $G$-space, there exists a $G$-map of the greatest ambit
$X$ to $Y$, and its restriction to $M$ is a $G$-map of $M$ to $Y$. Now let
$M'$ be another universal compact minimal $G$-space. There exist $G$-maps
$f:M\to M'$ and $g:M'\to M$. Since $M'$ is minimal, $f$ is surjective. On
the other hand, in virtue of Proposition~\ref{pr3} the composition
$gf:M\to M$ is bijective. It follows that $f$ is injective and hence a
$G$-isomorphism between $M$ and $M'$.
\end{proof}

Thus we have associated with every topological group $G$ the compact
$G$-space $M_G$. The question arises: what is this space? If $G$ is
discrete, then $M(G)$, being a retract of $S(G)=\beta G$, is extremally
disconnected. If $G$ is locally compact, the action of $G$ on $S(G)$ is
free \cite{Veech} (see also \cite[Theorem 3.1.1]{P3}), that is, if
$g\ne1$, then $gx\ne x$ for every $x\in S(G)$. It follows that the action
of $G$ on $M_G$ also is free. In some cases, the space $M(G)$ can be
described explicitly. For example, let $E$ be a countable infinite dicrete
space, and let $G=\Symm(E)\sbs E^E$ be the topologicall group of all
permutations of $E$. Then $M(G)$ can be identified with the space of all
linear orders on $E$ (Glasner--Weiss). Every linear order is considered as
a subset of $E\times E$, and the set of all subsets of $E\times E$ is
identified with the compact space $2^{E\times E}$. It is not clear whether
a similar result holds true if $E$ is uncountable.

Another example, due to V.Pestov, of a group $G$ for which $M(G)$ has an
explicit description is the following. Let $S^1$ be a circle, and let
$G=H_+(S^1)$ be the group of all orientation-preserving
self-homeomorphisms of $S^1$. Then $M_G$ can be identified with $S^1$
\cite[Theorem 6.6]{P1}. Pestov asked whether a similar assertion holds
for the Hilbert cube $Q=I^\omega$, where $I=[0,1]$: if $G=H(Q)$, are $M_G$
and $Q$ isomorphic as $G$-spaces?

The answer is no \cite{usptopproc}: there exists a compact $G$-space
$\Phi$ such that there is no $G$-map $Q\to \Phi$, hence $M(G)$ is not
isomorphic to $Q$ as a $G$-space. (I do not know whether $M(G)$ is
homeomorphic to $Q$, or whether $M(G)$ is metrizable.) One can take for
$\Phi$ the space of all maximal chains of closed subsets of $Q$. Thus
$\Phi\sbs \Exp\Exp Q$, where $\Exp K$ denotes the space of closed subsets
of a compact space $K$, equipped with the Vietoris topology.

A similar argument works in a more general situation. Let us say that the
action of a group $G$ on a $G$-space $X$ is {\em 3-transitive} if
$|X|\ge3$ and for any triples $(a_1,a_2,a_3)$ and $(b_1,b_2,b_3)$ of
distinct points in $X$ there exists $g\in G$ such that $ga_i=b_i$,
$i=1,2,3$. Suppose $K$ is compact and $G\sbs H(K)$ is a 3-transitive
group. Then there is no $G$-map from $K$ to $\Phi(K)$, the space of all
maximal chains of closed subsets of $K$. It follows that $M(G)\ne K$. This
argument implies

\begin{thm}[\cite{usptopproc}]
For every topological group $G$ the action of $G$ on the universal minimal 
compact $G$-space $M_G$ is not 3-transitive.
\label{main}
\end{thm}

For example, if $K$ is a compact manifold of dimension $>1$ or a compact
Menger manifold and $G=H(K)$, then $M(G)\ne K$, since the action of $G$ on
$K$ is 3-transitive. It would be interesting to understand what is $M(G)$
in this case.

Let $P$ be the pseudoarc (= the unique hereditarily indecomposable
chainable continuum) and $G=H(P)$. The action of $G$ on $P$ is transitive
but not 2-transitive, and the following question remains open:

\begin{Question}
Let $P$ be the pseudoarc and $G=H(P)$. Can $M_G$ be identified with~$P$?
\label{pseudo}
\end{Question}

Note that the argument involving the space $\Phi(K)$ of maximal chains
used above to prove that $M_G\ne K$ for every 3-transitive group $G\sbs
H(K)$ supports the conjecture that $M_G=P$ for $G=H(P)$, where $P$ is the
pseudoarc: there exists a $G$-map $P\to \Phi(P)$. Indeed, for every $x\in
P$ let $C_x$ be the collection of all subcontinua $F\sbs P$ such that
$x\in F$. Since any two subcontinua of $P$ are either disjoint or
comparable, $C_x$ is a chain. The chain $C_x$ can be shown to be maximal,
and the map $x\mapsto C_x$ from $P$ to $\Phi(P)$ is a $G$-map.

Pestov's example ($G=H_+(S^1)$, $M_G=S^1$) shows that the action of $G$ on
$M_G$ can be 2-transitive. Observe that there are precisely two $G$-maps
$S^1\to\Phi(S^1)$, which assign to every $x\in S^1$ the chain of all
closed arcs which either ``start at $x$'' or ``end at $x$'',
respectively.

The space $M_G$ is a singleton for many naturally arising non-locally
compact groups $G$. This property of $G$ is equivalent to the following
{\em fixed point on compacta (f.p.c.) property}: every compact $G$-space
has a $G$-fixed point. (A point $x$ is $G$-fixed if $gx=x$ for all $g\in
G$.) For example, if $H$ is a Hilbert space, the group $U(H)$ of all
unitary operators on $H$, equipped with the pointwise convergence
topology, has the f.p.c.\ property (Gromov-Milman); another example of a
group with this property, due to Pestov, is $H_+(\R)$, the group of all
orientation-preserving self-homeomorphisms of the real line. We refer the
reader to V.~Pestov's papers \cite{GP,P1,P2,P3,P4} on this subject.

Groups with the f.p.c.~property are also called {\em extremely amenable}.
Recall that a group $G$ is {\em amenable} if every continuous action of
$G$ by affine transformations on a convex compact subset of a locally
convex vector space has a $G$-fixed point. (This definition is equivalent
to the usual definition of amenability involving the existence of
invariant means). While every abelian topological group is amenable, it
may be or may be not extremely amenable. For example, a discrete group
$G\ne\{1\}$ is not extremely amenable; on the other hand, there exist
extremely amenable group topologies on the group $\Z$ of integers
\cite{Glas}. A necessary condition for a group $G$ to be extremely
amenable is that there be no non-constant continuous characters
$\chi:G\to\T$, where $\T=\{z\in\C:|z|=1\}$ is the unit circle. Indeed, if
$\chi:G\to\T$ is a character, $\chi\ne1$, then $G$ admits a fixed-point
free action on $\T$ given by $(g,x)\mapsto \chi(g)x$. It is not known
whether for abelian groups (or for the group $\Z$) the above necessary
condition is also sufficient:

\begin{Question}[Glasner]
Let $G$ be an abelian topological group. Suppose that $G$ has no 
non-trivial continuous characters $\chi:G\to \T$. 
Is $G$ extremely amenable?
\end{Question}

For cyclic groups the question can be reformulated as follows. Let $K$ be
a compact space, and let $f\in H(K)$ be a fixed-point free homeomorphism
of $K$. Let $G$ be the cyclic subgroup of $H(K)$ generated by $f$. Does
there exist a complex number $a$ such that $|a|=1$, $a\ne1$, and the
homomorphism $\chi:G\to \T$ defined by $\chi(f^n)=a^n$ is continuous?

If $K$ is a circle, the answer is yes: for every orientation-preserving
homeomorphism $f$ of a circle the rotation number is defined which gives
rise to a non-trivial continuous character on the group generated by $f$.

A positive answer to Glasner's question would imply the solution of the
following long-standing problem \cite{Glas, P2, P3}: is it true that for
every big set $S$ of integers the set $S-S$ contains a neighbourhood of
zero for the Bohr topology on $\Z$? A set $S$ of integers is said to be
{\em big} (or {\em syndetic}) if $S+F=\Z$ for some finite $F\sbs \Z$; this
means that the gaps between consequtive terms of $S$ are uniformly
bounded. The {\em Bohr topology} on $\Z$ is generated by all characters
$\chi:\Z\to\T$. It is known that for every big subset $S\sbs\Z$ the set
$S-S+S$ contains a Bohr neighbourhood of zero \cite[Corollary 3.25]{P3}.

Extremely amenable groups can be characterized in terms of big sets
\cite{P2, P3}. A subset $S$ of a topological group $G$ is {\em big on the
left}, or {\em left syndetic}, if $FS=G$ for some finite $F\sbs G$.

\begin{thm}[Pestov {\cite[Theorem 8.1]{P2}}]
A topological group $G$ is extremely amenable if and only if whenever
$S\sbs G$ is big on the left, $SS\obr$ is dense in $G$.
\end{thm}

There are other useful characterization of extremely amenable groups, also
due to Pestov \cite{P4, P5}, based on the notions of concentration and the
Ramsey-Dvoretzky-Milman property. The reader is invited to consult
\cite{P4}, \cite{P5} for details. We confine ourselves by formulating a
criterion of extreme amenability from \cite{P4}:

\begin{thm}[{\cite[Theorem 5.5]{P4}}]
A topological group $G$ is extremely amenable if and only if for every
bounded left uniformly continuous function $f$ from $G$ to a
finite-dimensional Euclidean space, every $\e>0$, and every finite (or
compact) $K\sbs G$ there exists $g\in G$ such that $\diam f(gK)<\e$.
\end{thm}

One of the main results in \cite{P4} is the following: the group $\Iso(U)$
is extremely amenable. Here $U$ is the Urysohn universal metric space. We
shall consider the group $\Iso(U)$ in Section~\ref{secury}.

\section{Roelcke compactifications}\label{secroel}

For a topological group $G$ let $R(G)$ be the maximal ideal space of the
$C^*$-algebra of all bounded complex functions on $G$ which are both left
and right uniformly continuous. The space $R(G)$ is the Samuel
compactification of the uniform space $(G, \sL\wedge\sR)$, where $\sL$ is
the left uniformity on $G$, $\sR$ is the right uniformity, and
$\sL\wedge\sR$ is the {\em Roelcke uniformity} on $G$, the greatest lower
bound of $\sL$ and $\sR$. We call $R(G)$ the {\em Roelcke
compactification} of $G$.

While the greatest lower bound of two compatible uniformities on a
topological space in general need not be compatible, the Roelcke
uniformity is compatible with the topology of $G$. The covers of the form
$\{UxU:x\in G\}$, $U\in \sN(G)$ constitute a base of uniform covers for
the Roelcke uniformity.

If $G$ is abelian, $R(G)=S(G)$. In general, $R(G)$ is a $G$-space, and the
identity map of $G$ extends to a $G$-map $S(G)\to R(G)$. A group $G$ is
{\em precompact} if one of the following equivalent properties holds:
$(G,\sL)$ is precompact; $(G,\sR)$ is precompact; $G$ is a subgroup of a
compact group. It can be shown that $G$ is precompact if and only if for
every neighbourhood $U$ of unity there exists a finite $F\sbs G$ such that
$G=FUF$. Let us say that $G$ is {\em Roelcke precompact} if the Roelcke
uniformity $\sL\wedge\sR$ is precompact. This means that for every
neighbourhood $U$ of unity there exists a finite $F\sbs G$ such that
$G=UFU$. There are many non-abelian non-precompact groups which are
Roelcke precompact. For example, the symmetric group $\Symm(E)$ of all
permutations of a discrete space $E$, or the unitary group $U(H)_s$ on a
Hilbert space $H$, equipped with the strong operator topology, are Roelcke
precompact. The Roelcke compactifications of these groups can be
explicitly described with the aid of the following construction.

Suppose that $G$ acts on a compact space $K$. For $g\in G$ let $\G(g)\sbs
K^2$ be the graph of the $g$-shift $x\mapsto gx$. The map $g\mapsto \G(g)$
from $G$ to $\Exp K^2$ is both left and right uniformly continuous (if the
compact space $\Exp K^2$ is equipped with its unique compatible
uniformity), hence it extends to a map $f_K:R(G)\to \Exp K^2$. If the
action of $G$ on $K$ is topologically faithful, the map $f_K$ often
happens to be an embedding, in which case $R(G)$ can be identified with
the closure of the set $\{\G(g):g\in G\}$ in $\Exp K^2$. For example, this
is the case if $K=S(G)$ or $K=R(G)$.

The space $\Exp K^2$ is the space of all closed relations on $K$. It has a
rich structure, since relations can be composed, reversed, or compared by
inclusion. This structure is partly inherited by $R(G)$. Let us consider
some examples.

\begin{example}
Let $G=\Symm(E)$ be the topological symmetric group. It acts on the
compact cube $K=2^E$. The natural map $f_K:R(G)\to \Exp K^2$ is an
embedding.
\end{example}

\begin{example}[{\cite{UHil}}]
Let $G$ be the unitary group $U(H)_s$ of a Hilbert space $H$, equipped
with the strong operator topology (this is the topology of pointwise
convergence inherited from the product $H^H$). Let $K$ be unit ball of
$H$. Equip $K$ with the weak topology. Then $K$ is compact. The unitary
group $G$ acts on $K$, and the map $R(G)\to \Exp K^2$ is an embedding. 

The space $R(G)$ has a better description in this case: $R(G)$ can be
identified with the unit ball $\Theta$ in the Banach algebra $\sB(H)$ of
all bounded linear operators on $H$. The topology on $\Theta$ is the weak
operator topology: the map $A\mapsto A|_K$ which assigns to every operator
of norm $\le1$ its restriction to $K$ is a homeomorphic embedding of
$\Theta$ into the compact space $K^K$. Thus $R(G)$ has a natural structure
of a semitopological semigroup. This can be used to deduce Stoyanov's
theorem: the group $G$ is minimal. Let us sketch the idea (see \cite{UHil}
for details). Let $f:G\to H$ be continuous homomorphism of $G$ onto a
topological group $H$. We want to prove that $f$ is open. To this end,
extend $f$ over $R(G)$. We get a map $F:\Theta\to R(H)$. Let
$S=F\obr(e_H)$ be its kernel. Then $S$ is a compact semigroup of
operators. If $S\sbs G$, then $G=F\obr(H)$ and $f$ is perfect, hence
quotient. For group homomorphisms `quotient' is equivalent to `open'. If
$S$ contains non-invertible operators, then $S$ contains idempotents (=
orthogonal projectors) other than 1. Since $S$ is invariant under inner
automorphisms of $\Theta$, it follows that $S$ contains 0, and this yields
$H=\{e_H\}$. 
\end{example}

\begin{example}[\cite{UHomeo}]
Let $K$ be a zero-dimensional compact space such that all non-empty clopen
subsets of $K$ are homeomorphic to $K$. (For example, $K$ may be the cube
$2^\kappa$ for some cardinal $\kappa$.) Let $G=H(K)$. The natural map
$f_K:R(G)\to \Exp K^2$ is an embedding. Moreover, the image of $f_K$,
which is the closure of the set of all graphs of self-homeomorphisms of
$K$, is the set $\Theta$ of all closed relations on $K$ whose domain and
range are equal to $K$. Thus $R(G)$ can be identified with $\Theta$.

This time $R(G)$ is an ordered semigroup, but not a semitopological
semigroup, since the composition of relations is not a separately
continuous operation. As in the previous example, one can use the space
$R(G)$ to prove that $G$ is minimal. Moreover, every non-constant onto
group homomorphism $f:G\to H$ is an isomorphism of topological groups. To
prove this, we proceed as before: extend $f$ to $F:\Theta\to R(H)$ and
look at the kernel $S=F\obr(e_H)$. Zorn's lemma implies the existence of
maximal idempotents in $S$ (with respect to inclusion). Symmetric
idempotents above the unity 1 (= the identity relation = the diagonal in
$K^2$) in $\Theta$ are precisely closed equivalence relations on $K$.
Since there are no non-trivial $G$-invariant closed equivalence relations
on $K$, there are no non-trivial choices for $S$: either $S=\{1\}$ or
$S=\Theta$. See \cite{UHomeo} for details. 
\end{example}

\begin{example}
Let $G=H_+(I)$ be the group of all orientation-preserving homeomorphisms
of the closed interval $I=[0,1]$. The map $f_G:R(G)\to \Exp I^2$ is a
homeomorphic embedding. Thus $R(G)$ can be identified with the closure of
the set of all graphs of strictly increasing functions $h:I\to I$ such
that $h(0)=0$ and $h(1)=1$. This closure consists of all curves $C\sbs
I^2$ which lead from $(0,0)$ to $(1,1)$ and look like graphs of increasing
functions, with the exception that $C$ may include both horizontal and
vertical segments.

There seems to be no natural semigroup structure on $R(G)$. This
observation leads to an important result, due to M.~Megrelishvili: the
group $G$ has no non-trivial homomorphisms to compact semitopological
semigroups and has no non-trivial representations by isometries in
reflexive Banach spaces. We discuss this result in the next section.
\end{example}

\begin{Question}
Let $G=H(Q)$, where $Q=I^\o$ is the Hilbert cube. Is the map
$f_Q: R(G)\to \Exp Q^2$ a homeomorphic embedding? Is the group $G$
minimal?
\end{Question}

\section{WAP compactifications}\label{secwap}

Let $S$ be a semigroup and a topological space. If the multiplication
$(x,y)\mapsto xy$ is separately continuous (this means that the maps
$x\mapsto ax$ and $x\mapsto xa$ are continuous for every $a\in S$), we say
that $S$ is a {\em semitopological semigroup}.

For a topological group $G$ let $f:G\to W(G)$ be the universal object in
the category of continuous semigroup homomorphisms of $G$ to compact
semitopological semigroups. In other words, $W(G)$ is a compact
semitopological semigroup, and for every continuous homomorphism $g:G\to
S$ to a compact semitopological semigroup $S$ there exists a unique
homomorphism $h:W(G)\to S$ such that $g=hf$.

The existence of $W(G)$ follows from two facts \cite[Ch.4]{Bour}: (1)
arbitrary products are defined in the category of compact semitopological
semigroups; (2) the cardinality of a compact space has an upper bound in
terms of its density. The space $W(G)$ can also be defined in terms of
weakly almost periodic functions. Recall the definition of such functions.

Let a topological group $G$ act on a space $X$. Denote by $C^b(X)$ the
Banach space of all bounded complex-valued continuous functions on $X$
equipped with the supremum norm. A function $f\in C^b(X)$ is called {\it
weakly almost periodic} ({\it w.a.p.} for short) if the $G$-orbit of $f$
is weakly relatively compact in the Banach space $C^b(X)$.

In particular, considering the left and right actions of a group $G$ on
itself, we can define left and right weakly almost periodic functions on
$G$. These two notions are actually equivalent \cite[Corollary
1.12]{Burckel}, so we can simply speak about w.a.p.\ functions on a group
$G$. The space WAP of all w.a.p.\ functions on $G$ is a $C^\ast$-algebra,
and the maximal ideal space of this algebra can be identified with $W(G)$.
Thus the algebra WAP is isomorphic to the algebra $C(W(G))$ of continuous
functions on $W(G)$. Call $W(G)$ the {\it weakly almost periodic {\rm
(}w.a.p.{\rm )} compactification} of the topological group~$G$.

\begin{remark}
In this section, by a {\it compactification} of a topological space $X$
we mean a compact Hausdorff space $K$ together with a continuous map
$j:X\to K$ with a dense range. We do {\it not} require that $j$ be a
homeomorphic embedding.
\end{remark}

For every reflexive Banach space $X$ there is a compact semitopological
semigroup $\Theta(X)$ associated with $X$: the semigroup of all linear
operators $A:X\to X$ of norm $\le1$, equipped with the weak operator
topology. Recall that a Banach space $X$ is reflexive if and only if the
unit ball $B$ in $X$ is weakly compact. If $X$ is reflexive, $\Theta(X)$
is homeomorphic to a closed subspace of $B^B$ (where $B$ carries the weak
topology) and hence compact.

It turns out that every compact semitopological semigroup embeds into
$\Theta(X)$ for some reflexive $X$:

\begin{thm}[Shtern \cite{Sh}, Megrelishvili \cite{MeElena}] 
Every compact semitopological semigroup is isomorphic to a closed
subsemigroup of $\Theta(X)$ for some reflexive Banach space $X$.
\label{ShtMeg}
\end{thm}

The group of invertible elements of $\Theta(X)$ is the group $\Iso_w(X)$
of isometries of $X$, equipped with the weak operator topology. This
topology actually coincides with the strong operator topology:

\begin{thm}[Megrelishvili \cite{MeElena}]
For every reflexive Banach space $X$ the weak and strong operator
topologies on the group $\Iso(X)$ agree.
\label{Megrefl}
\end{thm}

In particular, the group of invertible elements of $\Theta(X)$ is a
topological group. The natural action of this group on $\Theta(X)$ is
(jointly) continuous. This can be easily deduced from the fact (which
follows from Theorem~\ref{Megrefl}) that the topological groups
$\Iso_s(X)=\Iso_w(X)$ and $\Iso_s(X^*)=\Iso_w(X^*)$ are canonically
isomorphic. In virtue of Theorem~\ref{ShtMeg}, similar assertions hold
true for every compact semitopological semigroup $S$: the group $G$ of
invertible elements of $S$ is a topological group, and the map
$(x,y)\mapsto xy$ is jointly continuous on $G\times S$ (this is the
so-called Ellis-Lawson joint continuity theorem \cite{La}). Thus $S$ is
a $G$-space.

It follows that for every topological group $G$ the compact
semitopological semigroup $W(G)$ is a $G$-space, hence there exists a
$G$-map $S(G)\to W(G)$ extending the canonical map $G\to W(G)$. In terms
of function algebras this means that every w.a.p.\ function on $G$ is
right
uniformly continuous. Since the algebra WAP is invariant under the
inversion on $G$, w.a.p.\ functions are also left uniformly continuous and
hence Roelcke uniformly continuous. It follows that there is a natural map
$R(G)\to W(G)$.

If $G=U(H)$ is the unitary group of a Hilbert space $H$, then
$R(G)=\Theta(H)$ is a compact semitopological semigroup, and therefore the
canonical map $R(G)\to W(G)$ is a homeomorphism. Thus $W(G)=\Theta(H)$.
The canonical map $S(G)\to W(G)$ is a homeomorphism if and only if $G$ is
precompact \cite{MPU}.

In virtue of Theorems~\ref{ShtMeg} and~\ref{Megrefl}, the following two
properties are equivalent for every topological group $G$: (1) the
canonical map $G\to W(G)$ is injective; (2) there exists a faithful
representation of $G$ by isometries of a reflexive Banach space.
Similarly, the canonical map $G\to W(G)$ is a homeomorphic embedding if
and only if $G$ is isomorphic to a topological subgroup of $\Iso(X)$ for
some reflexive Banach space $X$. Does every topological group have these
properties? This long-standing question recently has been answered in the
negative by Megrelishvili:

\begin{thm}[\cite{MeForum}]
Let $G=H_+(I)$ be the group of all orientation-preserving homeomorphisms 
of $I=[0,1]$. Then $W(G)$ is a singleton. Equivalently, every
w.a.p.\ function on $G$ is constant.
\end{thm}

The proof is based on the description of the Roelcke compactification
$R(G)$ given in the previous section. Recall that $R(G)$ can be identified
with the space of ``monotonic curves''. Since there is a canonical onto 
map $R(G)\to W(G)$, the semigroup $W(G)$ can be obtained as a quotient of
$R(G)$. Megrelishvili proved that the only way to obtain a semitopological
semigroup from $R(G)$ is to collapse the whole space to a point.

\begin{Question}[Megrelishvili]
Does there exist a non-trivial {\em abelian} topological group $G$ for
which $W(G)$ is a singleton?
\end{Question}

\section{The group $\Iso(U)$}\label{secury}

In this section we consider a particular example of a topological group:
the group $\Iso(U)$ of isometries of the Urysohn universal metric space
$U$.

Let us say that a metric space $M$ is {\em $\o$-homogeneous} if every
isometry between two finite subsets of $M$ extends to an isometry of $M$
onto itself. A metric space $M$ is {\em finitely injective} if it has the
following property: if $K$ is a finite metric space and $L\sbs K$, then
every isometric embedding $L\to M$ can be extended to an isometric
embedding $K\to M$. The Urysohn universal space $U$ is the unique (up to
an isometry) complete separable metric space with the following
properties: (1) $U$ contains an isometric copy of any separable metric
space; (2) $U$ is $\o$-homogeneous. Equivalently, $U$ is the unique
finitely-injective complete separable metric space. The uniqueness of $U$
is easy: given two separable finitely-injective spaces $U_1$ and $U_2$,
one can use the ``back-and-forth'' (or ``shuttle'') method to construct an
isometry between countable dense subsets of $U_1$ and $U_2$; if $U_1$ and
$U_2$ are complete, they are isometric themselves. The existence of $U$
was proved by Urysohn \cite{Ury}; an easier construction was found by
Kat\v etov \cite{Kat}. If a metric on the set of integers is chosen at
random, the completion of the resulting metric space will be isometric to
$U$ with probability~1 \cite{Versh}.

Let $G=\Iso(U)$. The group $G$ is a universal topological group with a
countable base: every topological group $H$ with a countable base is
isomorphic (as a topological group) to a subgroup of $G$ \cite{UspUry}.
The idea of the proof is first to embed $G$ into $\Iso(M)$ for some
separable metric $M$ and then to embed $M$ into $U$ in such a way that
every isometry of $M$ has a natural extension to an isometry of $U$. Let
us give some details.

Our construction is based on Kat\v etov's paper \cite{Kat}. 
Let $(X,d)$ be a metric space.
We say that a function $f:X\to \R_+$ is {\it Kat\v etov\/} if
$|f(x)-f(y)|\le d(x,y) \le f(x)+f(y)$ for all $x,y\in X$. 
A function $f$ is Kat\v etov if and only if there exists a metric space
$Y=X\cup \{p\}$ containing $X$ as a subspace such that $f(x)$ for every
$x\in X$ is equal to the distance between $x$ and $p$. Let $E(X)$ be the
set of all Kat\v etov functions on $X$, equipped with the sup-metric. 
If $Y$ is a non-empty subset of $X$ and $f\in E(Y)$, define $g=\k_Y(f)\in
E(X)$ by
$$
g(x)=\inf\{d(x,y)+f(y):y\in Y\}
$$
for every $x\in X$.
It is easy to check that $g$ is indeed a Kat\v etov function on $X$ and
that $g$ extends $f$.
The map $\k_Y:E(Y)\to E(X)$ is an isometric embedding.
Let
$$
X^* = 
\bigcup\{\k_Y(E(Y)): Y\sbs X,\ Y\text{ is finite and non-empty}\,\}
\sbs E(X).
$$
For every $x\in X$ let $h_x\in E(X)$ be the function on $X$ defined by
$h_x(y)=d(x,y)$. Note that $h_x=\k_{\{x\}}(0)$ and hence $h_x\in X^*$.
The map $x\mapsto h_x$ is an isometric embedding of $X$ into $X^*$.
Thus we can identify $X$ with a subspace of $X^*$. 
If $K$ is a finite metric space, $L\sbs K$ and $|K\stm L|=1$, then every
isometric embedding of $L$ into $X$ can be extended to an isometric
embedding of $K$ into $X^*$.

Every isometry of $X$ has a canonical extension to an isometry of $X^*$, 
and we get an embedding of topological groups $\Iso(X)\to \Iso(X^*)$.
(Note that the natural homomorphism $\Iso(X)\to \Iso(E(X))$ in general 
need not be continuous.) Iterating the construction of $X^*$, we get an
increasing sequence of metric spaces $X\sbs X^*\sbs X^{**}\dots$.
Let $Y$ be the union of this sequence, and let $\cl{Y}$ be the completion 
of $Y$. 
We have a sequence of embeddings of topological groups
$$
\Iso(X)\to \Iso(X^*)\to \Iso(X^{**})\to\dots\to \Iso(Y)\to \Iso(\cl{Y}).
$$
The space $Y$ is finitely-injective.
The completion of a finitely-injective space is finitely-injective 
(Urysohn \cite{Ury}, see also \cite{minim}). 
Assume that $X$ is separable. 
Then $Y$ is separable, and $\cl{Y}$ is a complete separable
finitely-injective metric space.
Thus $\cl{Y}$ is isometric to $U$, and hence $\Iso(X)$ is isomorphic to a
topological subgroup of $\Iso(U)$.

Every topological group $G$ with a countable base is isomorphic to a
subgroup of $\Iso(X)$ for some separable Banach space $X$: there is a
countable subset $A\sbs\RUC(G)$ which generates the topology of $G$, and
we can take for $X$ the closed $G$-invariant linear subspace of $\RUC(G)$
generated by $A$. We just saw that $\Iso(X)$ is isomorphic to a subgroup
of $\Iso(U)$. Thus we have proved:

\begin{thm}[\cite{UspUry}]
Every topological group with a countable base is isomorphic to a
topological subgroup of the group $\Iso(U)$.
\end{thm}

Note that the group $\Iso(U)$ is Polish (= separable completely
metrizable). Another example of a universal Polish group is the group
$H(Q)$ of all homeomorphisms of the Hilbert cube \cite{UspFA}. To prove
that every topological group $G$ with a countable base is isomorphic to a
subgroup of $H(Q)$, it suffices to observe that: (1) $G$ is isomorphic to
a subgroup of $H(K)$ for some metrizable compact space $K$ (proof: we saw
that $G$ is isomorphic to a subgroup of $\Iso(B)$ for some separable
Banach space $B$, and we can take for $K$ the unit ball of the dual space
$B^*$, equipped with the $w^*$-topology); (2) if $K$ is compact and $P(K)$
is the compact space of all probability measures on $K$, there is a
natural embedding of topological groups $H(K)\to H(P(K))$; (3) if $K$ is
an infinite separable metrizable compact space, then $P(K)$ is
homeomorphic to the Hilbert cube. Until recently it remained unknown
whether the groups $H(Q)$ and $\Iso(U)$ are isomorphic or not. V.~Pestov
recently has proved that the group $\Iso(U)$ is extremely amenable
\cite{P4}. It follows that the groups $\Iso(U)$ and $H(Q)$ are {\em not}
isomorphic: the group $H(Q)$ is not extremely amenable, since the natural
action of $H(Q)$ on $Q$ has no fixed points.

A.M.~Vershik asked whether any two isomorphic compact subgroups of
$\Iso(U)$ are conjugate. The answer is in the negative: there exist
involutions $f,g\in\Iso(U)$ which are not conjugate. Moreover, $f$ has a
fixed point in $U$ and $g$ has no fixed points. The proof will appear
elsewhere.

The group $\Iso(U)$ is not Roelcke-precompact. To see this, fix $a\in U$
and consider the function $g\mapsto d(a,g(a))$ from $\Iso(U)$ to $\R_+$,
where $d$ is the metric on $U$. This function is $\sL\wedge\sR$-uniformly
continuous and unbounded, hence the Roelcke uniformity $\sL\wedge\sR$ is
not precompact. We slightly modify the space $U$, in order to obtain a
Roelcke-precompact group of isometries.

Let $U_1$ be the ``Urysohn universal metric space in the class of spaces
of diameter $\le1$''. This space is characterized by the following
properties: $U_1$ is a complete separable $\o$-homogeneous metric space of
diameter~1, and every separable metric space of diameter $\le1$ is
isometric to a subspace of $U_1$. Let $G=\Iso(U_1)$. This is a universal
Polish group. This group is Roelcke-precompact. Let us describe the
Roelcke compactification $R(G)$ of $G$.

Consider the compact space $K\sbs I^{U_1}$ of all non-expanding functions
$f:U_1\to I=[0,1]$. Then $K$ is a $G$-space, so there a natural map from
$R(G)$ to the set $\Exp K^2$ of all closed relations on $K$ (see
Section~\ref{secroel}). It turns out that this map is a homeomorphic
embedding.

There is a more geometric description of $R(G)$: it is the space of all
metric spaces $M$ of diameter~1 which are covered by two isometric copies
of $U_1$. More precisely, consider all triples $s=(M, i, j)$, where $M$ is
a metric space of diameter~1, $i:U_1\to M$ and $j:U_1\to M$ are isometric
embeddings, and $M=i(U_1)\cup j(U_1)$. Every such triple $s$ gives rise to
the function $p_s:U_1\times U_1\to I$ defined by $p_s(x,y)=d(i(x), j(y))$,
where $d$ is the metric on $M$. The set $\Theta$ of all functions $p_s$
that arise in this way is a compact subspace of $I^{U_1^2}$, and $R(G)$
can be identified with $\Theta$ \cite{minim}. Elements of $G$ correspond
to triples $(M,i,j)$ such that $M=i(U_1)=j(U_1)$.

The space $R(G)$ has a natural structure of an ordered semigroup (but not 
of a semitopological semigroup; it is likely that the w.a.p.\ 
compactification $W(G)$ of $G$ is a singleton).
If $R(G)$ is identified with a subset of $\Exp K^2$, then $R(G)$ happens
to be closed under composition of relations, whence the semigroup
structure, and the order is just the inclusion. 
If $R(G)$ is identified with $\Theta$, then the order is again natural,
and the semigroup operation is defined as follows: if $p,q\in \Theta$,
the product of $p$ and $q$ in $\Theta$ is the function $r:U_1^2\to I$
defined by
$$
r(x,y)=\inf(\{p(x,z)+q(z,y):z\in U_1\}\cup\{1\}) \quad (x,y\in U_1).
$$
There is a one-to-one correspondence between idempotents in $R(G)$ and 
closed subsets of $U_1$. The methods outlined in Section~\ref{secroel} 
imply the following:

\begin{thm}[\cite{minim}]
The universal Polish group $\Iso(U_1)$ is minimal.
\end{thm}

Thus every topological group with a countable base is isomorphic to a
subgroup of a minimal Roelcke-precompact Polish group. More generally,
every topological group is isomorphic to a subgroup of a minimal group of
the same weight \cite{minim}. The proof uses non-separable analogues of
the space $U_1$. Every topological group $G$ is isomorphic to a subgroup
of $\Iso(X)$, where $X$ is a complete $\o$-homogeneous metric space of
diameter~1 which is injective with respect to finite metric spaces of
diameter~1, and for every such $X$ the group $\Iso(X)$ is
Roelcke-precompact and minimal. The uniqueness of $X$ is lost in the
non-separable case, and it is not known whether there exists a universal
topological group of a given uncountable weight.

\providecommand{\bysame}{\leavevmode\hbox to3em{\hrulefill}\thinspace}
\providecommand{\MR}{\relax\ifhmode\unskip\space\fi MR }
\providecommand{\MRhref}[2]{%
  \href{http://www.ams.org/mathscinet-getitem?mr=#1}{#2}
}
\providecommand{\href}[2]{#2}

\end{document}